\documentclass[11pt]{amsart}
\usepackage{geometry}                % See geometry.pdf to learn the layout options. There are lots.
\usepackage[parfill]{parskip}    % Activate to begin paragraphs with an empty line rather than an indent
\usepackage{graphicx}
\usepackage{color}
\usepackage{amsmath,amscd}
\usepackage{amssymb}
\usepackage{epstopdf}
\usepackage[all]{xy}
\usepackage{eucal}
\usepackage{euler}
\usepackage{times}
\title{Pairings, duality, amenability and bounded cohomology}
 \def \note#1{{ \color{red}#1}}
\newtheorem{thm}{Theorem}
\newtheorem*{thm*}{Theorem}

\newtheorem*{mainthmthree}{Main Theorem}
\newtheorem*{BlockWeinbergerTheorem}{Theorem, Block and Weinberger}

\newtheorem{defn}{Definition}

\newtheorem{example}[defn]{Example}
\newtheorem{prop}[thm]{Proposition}

\newcommand{\C}{\mathbb{C}}

\newcommand{\imp}{$\implies$}
\newcommand{\norm}[1]{{\|{#1}\|}}

\newcommand{\skipit}[1]{\relax}

\def\note{\skipit}

\thanks{This research was partially supported by EPSRC grant  EP/F031947/1.}

\author{Jacek Brodzki}
\address{School of Mathematics, University of Southampton, Highfield, Southampton, SO17 1SH, England}
\email{J.Brodzki@soton.ac.uk}

\author{Graham A. Niblo}
\address{School of Mathematics, University of Southampton, Highfield, Southampton, SO17 1SH, England}
\email{G.A.Niblo@soton.ac.uk}

\author{Nick Wright}
\address{School of Mathematics, University of Southampton, Highfield, Southampton, SO17 1SH, England}
\email{N.J.Wright@soton.ac.uk}

\begin{document}
%%%%%%%%%%%%%%%%%%%%%%%%%%%%%%%%%%%%%%%%%%%%%%%%%%%%%%%%%%%%%%%%%%%%%%%%%%%%%
\begin{abstract}We give a new perspective on the homological characterisations of amenability given by Johnson in the context of bounded cohomology and by Block and Weinberger in the context of uniformly finite homology. We examine the interaction between their theories and explain the relationship between these characterisations. We apply these ideas to give a new proof of non-vanishing for the  bounded cohomology of a free group.\end{abstract}

\maketitle
%\section{}
%\subsection{}
%\begin{center}
%School of Mathematics, University of Southampton, Highfield, Southampton, SO17 1SH, England\\
%J.Brodzki@soton.ac.uk, G.A.Niblo@soton.ac.uk, N.J.Wright@soton.ac.uk.
%\end{center}

We will illuminate the relationship between the two following remarkable characterisations of amenability for a group. The definitions will follow the statements.

\begin{thm*} (Ringrose-Johnson, \cite{Johnson}) A group $G$ is amenable if and only if $H^1_b(G,(\ell^{\infty}(G)/\C)^*)=0$
\end{thm*}

\begin{thm*} (Block -Weinberger, \cite{BlockWeinberger}) A group $G$ is amenable if and only if $H_0^{uf}(G)\not = 0$.
\end{thm*}

It should be noted that both statements are part of a much larger picture. In the case of bounded cohomology vanishing of the first cohomology with the given coefficients is guaranteed by the triviality of a particular cocycle (the Johnson class $[J]\in  H^1(G, \ell^1_0G)$, defined below), and furthermore this ensures triviality of bounded cohomology with any coefficients in dimensions greater than or equal to $1$. In the case of Block and Weinberger's uniformly finite homology, vanishing of the zero dimensional homology group is guaranteed by the triviality of a fundamental class. It should also be noted that the the notion of amenability and the definition of uniformly finite homology can be extended from groups to arbitrary metric spaces, and the Block-Weinberger theorem applies in full generality. However there is no natural analog for the Ringrose and Johnson theorem in that context. We examine this issue further in \cite{BrodzkiNibloWright} where we define a cohomology (analogous to bounded cohomology) for a metric space. We use it there to give a generalisation of the Ringrose-Johnson theorem characterising Yu's property A, which is the natural generalisation of amenability in the context of coarse geometry.

Here we give a short, unified proof of the Johnson-Ringrose and Block-Weinberger theorems  by exploiting duality and the short exact sequence of coefficients
 
$$0\to \C \xrightarrow{\iota} \ell^\infty G \xrightarrow{\pi} \ell^\infty G/\C \to 0.$$

A crucial ingredient in the argument is that we can demonstrate the non-vanishing of a cohomology class by pairing it with a suitable homology class. 

It is well known that the free group of rank $2$ is non-amenable and therefore according to Ringrose and Johnson  the class $[J]$ is non-zero. As an application of the duality principle we demonstrate this non-vanishing by pairing $J$ with an explicit $\ell^1$-cycle. This construction should be compared with the argument in \cite{Mitsumatsu} for the linear independence of the Brooks cocycles in the bounded cohomology of a surface.

Recall the following definitions.
\begin{defn}
A \emph{mean} on a group $G$ is a positive linear functional $\mu$ on $\ell^\infty G$ such that $\norm{\mu}=1$. A group $G$ is \emph{amenable} if it admits a $G$-invariant mean.
\end{defn}

Recall that for a Banach space $V$ equipped with an isometric action of a group $G$, $C_b^m(G,V^*)$ denotes the $G$-module of equivariant bounded cochains $\phi:G^{m+1}\rightarrow V^*$. (Here bounded is defined by the Banach norm on the dual space $V^*$). This yields a cochain complex $(C_b^m(G,V^*), d)$ where $d$ denotes the natural differential induced by the homogeneous bar resolution. The cohomology of this complex is the bounded cohomology of the group with coefficients in $V^*$, denoted $H^*_b(G, V^*)$. For $V=\ell^\infty G/\C$ there is a particular class in dimension $1$ which detects amenability which we will call the Johnson element. This is represented by the function 

\[
J(g_0, g_1)=\delta_{g_1}-\delta_{g_0},
\]

where $\delta_g$ denotes the Dirac delta function supported at $g$. Note that $J(g_0, g_1)$ lies in the predual $\ell^1_0(G)$ of $V$, which we view as a subspace of its double dual, $V^*$.

Dually we have the chain complex $(C_m^{\ell^1}(G,V), \partial)$, where $C_m^{\ell^1}(G,V)$ consists of equivariant functions $c:G^{m+1}\rightarrow V$ which are $\ell^1$ on the subspace $\{e\}\times G^m$. The boundary map is defined by 

\[
\partial c(g_0, \ldots, g_{m-1})=\sum\limits_{g\in G, i\in\{0,\ldots, m\} } (-1)^ic(g_0, \ldots, g_{i-1}, g, g_i, \ldots, g_{m-1}).
\]
 
 The homology of this complex is the $\ell^1$-homology of the group with coefficients in  $V$, denoted $H^{\ell^1}_*(G, V)$. Note that there is a forgetful map $H_*(G, V)\rightarrow H^{\ell^1}_*(G, V)$. The fundamental class of Block and Weinberger in $H_0(G, \ell^\infty G)$ is represented by the cycle $c:G\rightarrow \ell^\infty G$ defined by $c(g)(h)=1$ for all $g,h\in G$. Applying the forgetful functor we obtain an element of $H^{\ell^1}_0(G, \ell^\infty G)$, and we will see that non-vanishing of this also characterises amenability. 
 
 We note that the pairing of $V^*$ with $V$, denoted $\langle-. -\rangle_V$ induces a pairing of $H^m_b(G, V^*)$ with $H_m^{\ell^1}(G, V)$ defined by 
 
 \[
 \langle[\phi], [c]\rangle = \sum\limits_{g_1,\ldots, g_m\in G}\langle\phi(e, g_1, \ldots, g_m), c(e, g_1, \ldots, g_m)\rangle_V.
 \]
 
It is clear that the pairing is defined at the level of cochains. To verify that it is well defined on classes one checks that the differential $d$ is the adjoint of the boundary map $\partial$.

The proof of the following result is a standard application of the snake lemma:

\begin{prop}\label{ses}
The short exact sequence of $G$-modules
$$0\to \C \xrightarrow{\iota} \ell^\infty G \xrightarrow{\pi} \ell^\infty G/\C \to 0.$$
induces a short exact sequence of chain complexes
$$0\to C_m^{\ell^1}(G,\C)\xrightarrow{\iota} C_m^{\ell^1}(G,\ell^\infty G) \xrightarrow{\pi} C_m^{\ell^1}(G,\ell^\infty G/\C)\to 0$$
and hence a long exact sequence of $\ell^1$-homology groups.

The short exact sequence of $G$-modules
$$0\to (\ell^\infty G/\C)^* \xrightarrow{\pi^*} \ell^\infty G^*  \xrightarrow{\iota^*} \C\to 0$$
induces a short exact sequence of cochain complexes
$$0\to C^m_{b}(G,(\ell^\infty G/\C)^*)\xrightarrow{\pi^*} C^m_{b}(G,\ell^\infty G^*) \xrightarrow{\iota^*} C^m_{b}(G,\C)\to 0$$
and hence a long exact sequence of bounded cohomology groups.
\qed
\end{prop}

Now we consider the Block-Weinberger uniformly finite homology of a countable discrete group $G$, where the group is equipped with a proper left invariant metric $d$. All of the definitions in \cite{BlockWeinberger} work for an arbitrary metric space, however we will restrict attention to the world of groups. Here the definition is slightly simplified as we have coarse bounded geometry and, as we shall see, it is natural to relate uniformly finite homology to bounded cohomology.

Let $C^{uf}_q(G, \mathbb R)$ denote the vector space of real valued functions $\phi:G^{q+1}\rightarrow \mathbb R$ which are bounded, and have controlled support. That is to say there is a constant $K$ (depending on the function $\phi$) such that if $\text{diam}\{g_0, \ldots, g_q\}\geq K$ then $\phi(g_0, \ldots, g_q)=0$.
The differential $\partial$ on the homogeneous bar resolution defined by 
\[
\partial(g_0, \ldots, g_q)= \sum\limits_{i=0}^q(-1)^i(g_0, \ldots, \hat{g_i}, \ldots, g_q)
\]
  extends linearly to induce a chain map $\partial: C^{uf}_q(G, \mathbb R)\rightarrow C^{uf}_{q-1}(G, \mathbb R)$. The uniformly finite homology of $G$ is then the homology of this chain complex. Block and Weinberger showed that there is a fundamental class $[\mathbf 1]$ in degree $0$ represented by the constant function $g\mapsto 1$ which detects amenability.

\begin{BlockWeinbergerTheorem}[\cite{BlockWeinberger} Theorem 3.1] The group $H^{uf}_0(G, \mathbb R)=0$ if and only if the fundamental class $[\mathbf 1]=0$, if and only if $G$ is not amenable.
\end{BlockWeinbergerTheorem}

In fact the uniformly finite homology coincides with the classical group homology $H_q(G, \ell^\infty(G))$, with coefficients in the module of bounded real valued functions on $G$. The corresponding chain complex consists of functions $\phi:G^{q+1}\rightarrow \ell^\infty(G)$ which are equivariant and supported on finitely many $G$-orbits. 
To see that the two homologies coincide we note that a cochain $\phi\in  C^{uf}_q(G, \mathbb R)$ can be inflated to a map $\overline \phi:G^q\rightarrow \ell^\infty(G)$ by setting $\overline \phi(g_0, \ldots, g_q)(g)=\phi(g^{-1}g_0, \ldots, g^{-1}g_q)$. This function is, by construction, equivariant and the controlled support condition  ensures that $\overline \phi$ is a chain in the  group  homology chain complex with $\ell^\infty$ coefficients. It is easy to see that this process is invertible. The differentials in both complexes are induced by the homogeneous bar resolution, so this map is an isomorphism between the chain complexes. Hence we may identify $H^{uf}_q(G, \mathbb R)$ with $H_q(G, \ell^\infty(G))$.

Let $\mathbf 1$ denote the constant function $G\rightarrow \mathbb C$ which takes the value $1$ at every $g\in G$. This function represents classes in all of the following objects: $H^0_b(G, \C)$, $H_0(G, \C)$, $H^{\ell^1}_0(G, \C)$. Our point of view is that the Block-Weinberger fundamental class is $i[\mathbf 1]\in H_0(G, \ell^\infty G)$, while the Johnson cocycle is $d[\mathbf 1]\in H^1_b(G, (\ell^\infty G/\C)^*)$, where $d$ denotes the connecting map $H^0_b(G, \C)\rightarrow H^1_b(G, (\ell^\infty G/\C)^*)$. The first of these observations is elementary. For the second, note that $d[\mathbf 1]$ is obtained by lifting $\mathbf 1$ to the element $g\mapsto \delta_g$ in $C^0_b(G,(\ell^\infty G)^*)$ and taking the coboundary. This produces the Johnson cocycle $J(g_0,g_1)=\delta_{g_1}-\delta_{g_0}$.

By exploiting the connecting maps arising in Proposition \ref{ses} together with these observations we will obtain a new proof that $G$ is amenable if and only if the Johnson cocycle in bounded cohomology vanishes, and that this is equivalent to non-vanishing of the Block-Weinberger fundamental class. The first hint of the interaction is provided by the duality between $H_0(G, \ell^\infty G)$ and $H^0(G, \ell^\infty G^*)$, and the observation that the latter is equal to $H^0_b(G, \ell^\infty G^*)$ since equivariance ensures that 0-cochains are bounded. The non-vanishing of $H^0_b(G, \ell^\infty G^*)$ is equivalent to amenability, since elements of  $H^0_{b}(G,\ell^\infty G^*)$ are maps $\phi:G\rightarrow \ell^\infty G^*$, which are $G$-equivariant and also, since they are cocycles, constant on $G$. Hence the value of a cocycle $\phi$ at any (and hence all) $g\in G$ is a $G$-invariant linear functional on $\ell^\infty G$. If $\phi$ is non-zero then taking its absolute value and normalising we obtain an invariant mean on the group. Conversely any invariant mean on the group is an invariant linear functional on $\ell^\infty G$ and hence gives a non-zero element of $H^0_{b}(G,\ell^\infty G^*)$.

\begin{mainthmthree}\label{beautiful}
Let $G$ be a countable discrete group. The following are equivalent:
\begin{enumerate}
\item \label{amen} $G$ is amenable.

\item \label{inv-mean} $\iota^*:H^0_{b}(G,\ell^\infty G^*) \to H^0_{b}(G,\C)$ is surjective.

\item \label{johnson} The Johnson class $d[\mathbf 1]$ vanishes in $H^1_{b}(G,(\ell^\infty G/\C)^*)$. %(Note: $d[\mathbf 1]$ is the Johnson cocycle.)

%\item \label{inv-funct} $H^0_{b}(G,\ell^\infty G^*)$ is non-zero.

\item \label{pairing} $\langle d[\mathbf 1],[c]\rangle =0$ for all $[c]$ in $H_1^{\ell^1}(G,\ell^\infty G/\C)$. (Hence for a non-amenable group, the non-triviality of $d[\mathbf 1]$ is detected by the pairing.)

\item \label{hom-last} $\iota[\mathbf 1] \in H_0^{\ell^1}(G,\ell^\infty G)$ is non-zero.

%\item \label{a-inv-funct} $d:C^0_{b}(G,\ell^\infty G^*)\to C^1_{b}(G,\ell^\infty G^*)$ is not bounded below.

\item \label{B-W}The Block-Weinberger fundamental class $\iota[\mathbf 1] \in H_0(G,\ell^\infty G)$ is non-zero. %(Note: $i[\mathbf 1]$ is the Block-Weinberger cycle.)

%\item \label{a-inv-funct} $d:C^0_{}(G,\ell^\infty G^*)\to C^1_{}(G,\ell^\infty G^*)$ is not bounded below.

\end{enumerate}
\end{mainthmthree}

\begin{proof}
(\ref{amen})\imp (\ref{inv-mean}) since $H^0_{b}(G,\C)=\C$, and for $\mu$ an invariant mean $i^*[\mu]=[\mathbf 1]$.

(\ref{inv-mean}) $\iff$ (\ref{johnson}): By exactness, surjectivity of $\iota^*$ is equivalent to vanishing of $d$, hence in particular this implies $d[\mathbf 1]=0$. The converse follows from the fact that $[\mathbf 1]$ generates $H^0_{b}(G,\C)$, so if $d[\mathbf 1]=0$ then $d=0$ and $\iota^*$ is surjective.

The implication %(\ref{inv-mean})\imp (\ref{inv-funct}) and 
(\ref{johnson})\imp (\ref{pairing}) is trivial.

%(\ref{inv-funct})\imp(\ref{a-inv-funct}) If $H^0_{b}(G,\ell^\infty G^*)$ is non-zero then $H^0_{}(G,\ell^\infty G^*)$ is non-zero because in dimension $0$ the forgetful functor is an isomorphism, from which it is immediate that $d:C^0_{}(G,\ell^\infty G^*)\to C^1_{}(G,\ell^\infty G^*)$ is not injective and therefore not bounded below.

(\ref{pairing}) $\implies$ (\ref{hom-last}): (\ref{pairing}) is equivalent to $\langle [\mathbf 1],\partial[c]\rangle =0$ for all $[c]$ in $H_1^{\ell^1}(G,\ell^\infty G/\C)$ by duality. We note that the space of 0-cycles in $C_0^{\ell^1}(G,\C)$ is $\C$, and noting that the pairing of the class $[\mathbf 1]$ in $H^0_{b}(G,\C)$ with the class $[\mathbf 1]$ in $H_0^{\ell^1}(G,\C)$ is $\langle [\mathbf 1],[\mathbf 1]\rangle =1$, we see that $[\mathbf 1]\in H_0^{\ell^1}(G,\C)$ is not a boundary. Thus $H_0^{\ell^1}(G,\C)=\C$ and the pairing with $H^0_{b}(G,\C)$ is faithful so $\langle [\mathbf 1],\partial[c]\rangle =0$ for all $[c]$ implies $\partial=0$. From this we deduce that $\iota$ is injective by exactness, hence we have (\ref{hom-last}): $\iota[\mathbf 1]$ is non-zero.

%(\ref{hom-last}) $\implies$ (\ref{a-inv-funct}):  (\ref{hom-last}) implies that $\partial:C_1^{\ell^1}(G,\ell^\infty G)\to C_0^{\ell^1}(G,\ell^\infty G)$ is not surjective, since in particular $i[\mathbf 1]$ is not in the image of $\partial$. This implies (\ref{a-inv-funct}) since if $d=\partial^*$ is bounded below then $d^*=\partial^{**}$ is surjective, whence $\partial$ is also surjective.

%(\ref{a-inv-funct})\imp(\ref{amen}) since (\ref{a-inv-funct}) is the existence of a sequence of linear functionals $\mu_n$ on $\ell^\infty G$ such that $\|\mu_n\|=1$ for all $n$ and $\|g\mu_n-\mu_n\|$ converges to zero uniformly in $g$ as $n\to\infty$, which implies amenability by Theorem \ref{a-inv-funct-thm}. (In fa

(\ref{hom-last})\imp (\ref{B-W}) since $\iota[\mathbf 1] \in H_0^{\ell^1}(G,\ell^\infty G)$ is the image of the corresponding element of $H_0(G,\ell^\infty G)$ under the forgetful map.

(\ref{B-W})\imp (\ref{amen}): We will use an argument due to Nowak. Let $\delta:C^0_{}(G,\ell^1(G))\to C^1_{}(G,\ell^1(G))$ denote the restriction of $d$. This is the predual of $\partial$. First we note that  $\delta$ is not bounded below, since if it were then  $\partial=\delta^*$ would be surjective and $H_0(G, \ell^\infty G)$ would vanish giving $\iota[\mathbf 1]=0$, which is a contradiction.

The fact that $\delta$ is not bounded below is precisely the assertion that there is a Reiter sequence for the group and that therefore it is amenable. \note{Do we need to say what a Reiter sequence is?}

\end{proof}

As an example of this approach we give a proof of non-amenability for $F_2$ by constructing an explicit element $[c]\in H_1^{\ell^1}(G,\ell^\infty G/\C)$ for which $\langle d[\mathbf 1],[c]\rangle\not= 0$.

Let $\{a,b\}$ be a free basis for $F_2$, and let $\Gamma$ denote the Cayley graph of $F$ with respect to this generating set. $\Gamma$ is a tree and the action of $G$ on $\Gamma$ extends to the Gromov boundary. We choose a point $p$ in the Gromov boundary of $\Gamma$. For the sake of definiteness we set $p$ to be the endpoint of the ray $(a^n)$ where $n$ ranges over the positive integers, though this is not essential.

For a generator $s$ of $F_2$ (or its inverse) we set $c(e,s)(g)= 1$ if $(e, s)$ is the first edge on the geodesic from $e$ to $gp$ and set $c(e, s)(g)=0$ otherwise. Extending the definition by equivariance we obtain a function $c$ defined on the edges of $\Gamma$ with values in $\ell^\infty G$ and this represents an element $\overline c\in\ell^\infty G/\mathbb C$. 

Now consider $\partial c(e)=\sum\limits_{s\in\{a^{\pm 1}, b^{\pm 1}\}}c(s,e)-c(e,s)$.

For a given $g$ exactly one of the edges $(e,a), (e,b), (e, a^{-1}), (e, b^{-1})$ is the first edge on the geodesic $[e, gp]$, so the sum  $c(e,a)+ c(e,b)+ c(e, a^{-1})+c(e, b^{-1})$ is the constant function $\mathbf 1$ on $G$. 

On the other hand for a generator $s$, $c(s,e)(g)=1$ if and only if the edge $(s, e)$ is the first edge on the geodesic from $s$ to $gp$. We now consider the function $c(a,e)+c(b,e)+c(a^{-1}, e)+c(b^{-1},e)$. For a given $g\in G$ there is a unique point in the set $\{a,b, a^{-1}, b^{-1}\}$ which lies on the geodesic from $e$ to $gp$, and this is the only one for which the corresponding term of the sum takes the value $0$, so the sum $c(a,e)+c(b,e)+c(a^{-1}, e)+c(b^{-1},e)$ is the constant function $\mathbf 3$.

Hence $\partial c(e)=\mathbf 3-\mathbf 1=\mathbf 2$. Now by equivariance $\partial c(k)=\mathbf 2$ for all $k$, hence $\partial \overline c$ vanishes in $\ell^\infty G/\mathbb C$, so $\overline c$ is  a cycle and therefore represents an element $[\overline c]\in H_1^{\ell^1}(G,\ell^\infty G/\C)$.

We now compute the pairing $\langle d[\mathbf 1], [\overline c]\rangle$. 

\[
\langle d[\mathbf 1], [\overline c]\rangle=\langle [\mathbf 1], \partial [\overline c]\rangle=\langle [\mathbf 1], [\partial  c]\rangle = \langle [\mathbf 1], [\mathbf 2]\rangle=2.
\]

%\begin{align*}
%\langle d[\mathbf 1], [c]\rangle &= \sum\limits_{g\in \{a^{\pm 1}, b^{\pm 1}\}}\sum\limits _{h\in G} c(e,g)(h)(\delta_g-\delta_e)(h)\\
%&=c(e,a)(a)-c(e,a)(e)+c(e,b)(b)-c(e,b)(e)\\
%&+c(e, a^{-1})(a^{-1})-c(e,a^{-1})(e)+c(e, b^{-1})(b^{-1})-c(e,b^{-1})(e)\\
%&=1-1+1-0+0-0+1-0
%=2\\
%\end{align*}

Hence $F_2$ is not amenable.

We conclude by noting that amenability is also equivalent to vanishing of the Johnson class as an element of the classical group cohomology  $H^1(G,(\ell^\infty G/\C)^*)$. To see this, replace the pairing of    $H^1_{b}(G,(\ell^\infty G/\C)^*)$  and $H_1^{\ell^1}(G,\ell^\infty G/\C)$ in the proof of Theorem \ref{beautiful} with the standard pairing of $H^1(G,(\ell^\infty G/\C)^*)$  and $H_1(G,\ell^\infty G/\C)$, hence deducing that vanishing of the Johnson element in $H^1(G,(\ell^\infty G/\C)^*)$ implies non-vanishing of the Block-Weinberger fundamental class. Hence we obtain the following theorem.

\begin{thm}\label{unboundedbeauty}
Let $G$ be a countable discrete group. The following are equivalent:
\begin{enumerate}
\item $G$ is amenable.

\item  $\mathbf 1$ lies in the image of $i^*:H^0_{}(G,\ell^\infty G^*) \to H^0_{}(G,\C)$.

\item The Johnson class $d[\mathbf 1]$ vanishes in $H^1_{}(G,(\ell^\infty G/\C)^*)$. %(Note: $d[\mathbf 1]$ is the Johnson cocycle.)

\end{enumerate}

\end{thm}


\begin{thebibliography}{999}

\bibitem{BlockWeinberger} J.~Block and S.~Weinberger, Aperiodic tilings, positive scalar curvature and amenability of spaces. J. Amer. Math. Soc. 5 (1992), no. 4, 907--918.
\bibitem{BrodzkiNibloWright} J.~Brodzki, G.~A.~Niblo and N.~J.~Wright, A cohomological characterisation of Yu's property A for metric spaces. arXiv:1002.5040 (2010).
\bibitem{Brooks} R.~Brooks, Some remarks on bounded cohomology. Riemann surfaces and related topics: Proceedings of the 1978 Stony Brook Conference (State Univ. New York, Stony Brook, N.Y., 1978), pp. 53--63, Ann. of Math. Stud., 97, Princeton Univ. Press, Princeton, N.J., (1981). 
\bibitem{Johnson} B.~E.~Johnson, Cohomology of Banach Algebras, Memoirs of the AMS Number 127, 1972, AMS, Providence, Rhode Island.
\bibitem{Mitsumatsu} Y.~Mitsumatsu, Bounded Cohomology and $\ell^1$ Homology of Surfaces, Topology Vol. 23. No. 4, (1984) pp. 465-471.
\bibitem{Yu} G.~Yu, The {C}oarse {B}aum-{C}onnes conjecture for spaces which admit a uniform embedding into {H}ilbert space, Inventiones Math. \textbf{139}
  (2000), 201--240.
  \end{thebibliography}
\end{document}